\pdfoutput=1

\documentclass[twocolumn]{autart}    

\usepackage{graphicx}          
\usepackage{mathematix}
\usepackage{hyphenat}
\usepackage{hyperref}
\usepackage{calc}
\usepackage[inline]{enumitem}
\newcommand{\ttop}{{\scriptstyle\top}}
\newcommand{\tttop}{{\scriptscriptstyle\top}}
\newcommand{\tplus}{{\scriptscriptstyle+}}

\newcommand{\fcns}{\operatorname{fcns}}
\usepackage{flushend}
\graphicspath{{Pics/}}

\newcommandx\Cexp[4][1={},3={},4={}]{%
	\ifstrempty{#3}{%
    	\E_{#1}\ifstrempty{#2}{}{\left[{}#2{}\ifstrempty{#4}{}{{};{}#4}\right]}
    }{%
    	\E_{#1}\ifstrempty{#2}{}{{\left[{}#2{}\right|}{\left.\vphantom{#2}#3\ifstrempty{#4}{}{{};{}#4}\right]}}
    }
}
\newcommandx\Crho[4][1={},3={},4={}]{%
    \ifstrempty{#3}{%
    	\rho_{#1}\ifstrempty{#2}{}{\left[{}#2{}\ifstrempty{#4}{}{;#4}\right]}
    }{%
    	\rho_{#1}\ifstrempty{#2}{}{{\left[{}#2{}\right|}{\left.\vphantom{#2}#3\ifstrempty{#4}{}{{};{}#4}\right]}}
    }
}

\usepackage{acronym}
\acrodef{mpc}[MPC]{model predictive control}
\acrodef{mjls}[MJLS]{Markov jump linear systems}
\acrodef{qcqp}[QCQP]{quadratically constrained quadratic program}
\acrodef{lmi}[LMI]{linear matrix inequality}
\acrodef{ui}[UI]{uniformly invariant}
\acrodef{dp}[DP]{dynamic programming}
\acrodef{rses}[RSES]{risk-square exponentially stable}
\acrodef{lsc}[lsc]{lower semi-continuous}

\makeatletter
\newcommand{\raisemath}[1]{\mathpalette{\raisem@th{#1}}}
\newcommand{\raisem@th}[3]{\raisebox{#1}{$#2#3$}}
\makeatother

\begin{document}

\begin{frontmatter}

\title{Risk-averse model predictive control\thanksref{footnoteinfo}}

\thanks[footnoteinfo]{This is a preprint; please, cite the published version of the 
paper as follows: P.~Sopasakis, D.~Herceg, A.~Bemporad and P.~Patrinos, Risk-averse model predictive control,
Automatica, 100:281--288, 2019. DOI: \texttt{10.1016/j.automatica.2018.11.022}.}

\author[KIOS]{Pantelis~Sopasakis}\ead{sopasakis.pantelis@ucy.ac.cy},%
\author[IMTL]{Domagoj~Herceg},%
\author[IMTL]{Alberto Bemporad},%
\author[KUL]{Panagiotis Patrinos}%

\address[KIOS]{University of Cyprus, 
               Department of Electrical and Computer Engineering,
               KIOS Research Center for Intelligent Systems and Networks, 
               1 Panepistimiou Avenue, 2109 Aglantzia, Nicosia, Cyprus.}
\address[IMTL]{IMT School for Advanced Studies Lucca, Piazza San Francesco 19, 55100 Lucca, Italy.}                                           
\address[KUL]{KU Leuven, Department of Electrical Engineering (ESAT), STADIUS Center for 
               Dynamical Systems, Signal Processing and Data Analytics, 
               Kasteelpark Arenberg 10, 3001 Leuven, Belgium.}


\begin{abstract}
Risk-averse model predictive control (MPC) offers a control framework that 
allows one to account for ambiguity in the knowledge of the underlying 
probability distribution and unifies stochastic and worst-case MPC.  
In this paper we study risk-averse MPC problems for constrained nonlinear 
Markovian switching systems using generic cost functions, and derive 
Lyapunov-type risk-averse stability conditions by leveraging the properties of 
risk-averse dynamic programming operators. We propose a controller design
procedure to design risk-averse stabilizing terminal conditions for constrained 
nonlinear Markovian switching systems. Lastly, we cast the resulting 
risk-averse optimal control problem in a favorable form which can be solved 
efficiently and thus deems risk-averse MPC suitable for applications.
\end{abstract}%

\end{frontmatter}
\section{Introduction}
\subsection{Background and Contributions}
There exist two main ways to deal with uncertainty in \ac{mpc}, 
namely, the \textit{robust} and the 
\textit{stochastic} approaches.
In \textit{robust \ac{mpc}}, modeling errors or disturbances are 
modeled as unknown-but-bounded quantities and the performance index 
is minimized with respect to the worst-case realization
(min-max approach)~\cite{rawlings2009model}. However, such worst-case 
events are unlikely to occur in practice and render robust \ac{mpc} 
severely conservative since all statistical information, typically 
available from past measurements, is ignored.

On the other hand, in \textit{stochastic \ac{mpc}} we assume that the underlying 
uncertainty is a random vector following some probability distribution~\cite{smpc-mesbah}
and minimize the expectation of a performance index; such formulations
are significantly less conservative.
The driving random process is often taken to be 
normally and independently identically distributed~\cite{Hokayem201277} or
it is assumed that it is a finite Markov process~\cite{smpc4mss} and in 
\textit{scenario-based \ac{mpc}}, filtered probability distributions are 
estimated from data~\cite{HanSopBemRaiCol15}. 
However, not always can we accurately estimate 
a distribution from available data, nor does it remain 
constant in time. Stochastic \ac{mpc} will guarantee mean-square stability 
of the closed-loop system only with respect to the nominal probability 
distribution, therefore, errors in the estimation of that distribution 
may lead to bad performance or even instability.

The theory of \textit{risk measures}~\cite{shapirolectures}
allows to interpolate between these two extreme cases.
Roughly speaking, risk measures quantify the importance and effect 
of the right tail of a distribution of losses, that is, the impact of the 
occurrence of \textit{extreme events}. As such they offer a 
mathematically elegant tool to tackle problems where we seek 
to avoid \textit{high effect low probability} (HELP) events and can be readily used 
in various applications.

The first steps to risk-averse formulations can be traced back 
to linear-exponential-quadratic Gaussian control~\cite{leqg73}
and the study of stochastic control problems under inexact knowledge of 
the underlying probability distribution which is often termed
\textit{distributionally robust}~\cite{Goh2010}. 
Distributionally robust control methodologies have been proposed 
for linear systems with probabilistic constraints
assuming knowledge of some moments of the distribution%
~\cite{PKGM:2016}. 
The same problem was also recently addressed for Markov decision 
processes with uncertain transition probabilities~\cite{YuXu16}.

Risk-averse \ac{mpc} formulations for \ac{mjls}
are studied in~\cite{chow2017arxiv,chow-risk-averseMPC}. In~\cite{chow-risk-averseMPC} 
the authors formulate an \ac{mpc} optimization problem employing a coherent risk measure of an
uncertain cost as an objective function and  give conditions under which 
the \ac{mpc} control law is stabilizing, albeit for a system with no state-input constraints. 
This is extended in~\cite{chow2017arxiv} assuming ellipsoidal state-input constraints. 
Building up on these results, we further improve on the state of the art by studying nonlinear systems and proposing 
a computationally favorable formulation for risk-averse optimization problems which leads to low computation times.

In the optimization and operations research communities, the solution 
of multistage risk-averse optimal control problems
has been considered prohibitive as only bundle and cutting-plane
methods are currently used~\cite{Asamov2015,Collado2012ScenarioDO,Bruno2016979}.
Reported results are limited to short prediction horizons and linear 
stage cost functions.
An alternative solution approach solves the \ac{dp} problem
using multiparametric piecewise quadratic programming~\cite{patrinos2011convex}, but its applicability 
is limited to systems with few states and small prediction horizons~\cite{patrinosECC2007}. 
In a 2017 paper, Rockafellar proposed an algorithmic scheme for solving multistage
problems using a non-composite (not nested) risk measure  
recognizing the difficulty of solving problems with nested risk mappings%
~\cite{Rockafellar2017}.
Indeed, the difficulty lies in that the cost function is written as a series of compositions of 
typically nonsmooth operators.
In Section~\ref{sec:computational} we present a computationally 
tractable approach for the solution of multistage risk-averse 
problems by disentangling this series of compositions.
This formulation renders risk-averse \ac{mpc} suitable for embedded applications.

In this paper we formulate multistage risk-averse optimal control problems 
using Markov risk measures in a \ac{dp} setting and
derive Lyapunov-type risk-averse stability conditions.
We study risk-averse \ac{mpc} formulations
for nonlinear Markovian switching systems under generally nonconvex
joint state-input constraints and propose a 
controller design procedure for nonlinear systems
with smooth dynamics and Lipschitz-continuous gradient.
Lastly, we provide simulation examples to demonstrate 
the applicability of the proposed approach.

\subsection{Notation}
Let \(\barre \dfn \Re\cup\{+\infty\}\) be the set of extended-real numbers,
\(\N_{[k_1,k_2]}\) the integers in \([k_1,k_2]\),
for \(z\in\Re^n\) let \(\plus{z}{=}\max\{0,z\}\) (where the max is taken element-wise). 
We denote by \(1_n\) the vector in \(\Re^n\) with all coordinates equal to \(1\).
We denote the sets of \(n\)-by-\(n\) symmetric positive definite (semidefinite) matrices
as \(\Spp^n\) (\(\Sp^n\)).
For two \(n\)-by-\(n\) symmetric matrices \(M_1,M_2\), \(M_1\succcurlyeq M_2\)  
means that \(M_1{-}M_2\in \Sp^n\). 
We denote the transpose of a matrix \(A\) by \(A^{\tttop}\) and the identity 
matrix by \(I\). For a \(g:\Re^n\to\Re^m\), its \textit{Jacobian matrix} is the mapping
\(Jg:\Re^n\to \Re^{m\times n}\) defined as 
\(Jg(x) = (\nicefrac{\partial g_i(x)}{\partial x_j})_{i,j}\),
provided that the partial derivatives exist.
For \(\epsilon{}\geq{}0\) we define \(\mathcal{B}_\epsilon=\{x {}\mid{} \|x\|\leq \epsilon\}\).
For a set \(C\subseteq \Re^n\), we define its \textit{indicator function} 
as \(\delta_C(x) = 0\) if \(x\in C\) and \(\delta_C(x)=\infty\) otherwise.
The \textit{domain} of an extended-real-valued function \(f{:}\Re^n{\to}\barre\)
is \(\dom f = \set{x\in \Re^n}[f(x) {<} \infty]\).
An extended-real-valued function \(f:\Re^n\to\barre\) is called
\textit{proper} if its domain is nonempty;
it is called \ac{lsc} if its lower level sets are closed.
An \(\ell:\Re^n\times\Re^m\ni(x,u)\mapsto \ell(x,u)\in\barre\) is called \textit{level bounded
in \(u\) locally uniformly in \(x\)} if for each \(x_0\in\Re^n\) and \(\alpha\in\Re\), there is
a neighborhood \(U_{x_0}\) of \(x_0\) along with a bounded set \(B\subseteq \Re^m\) such that
\(\set{u}[\ell(x,u)\leq \alpha] \subseteq B\) for all \(x_0\in U_{x_0}\).
The \textit{effective domain} of a set-valued mapping \(F:\Re^n\rightrightarrows \Re^m\)
is defined as \(\dom F = \set{x\in\Re^n}[F(x) \neq \emptyset]\).
For a nonempty set \(E\) and a finite set \(\NN\) we define \(\fcns(E,\NN) = \set{V:E \times \NN \to \barre}	
[V(x,i) \geq 0, V(0,i)=0, \text{ for all } x\in E, i \in \NN]\).

\section{Risk-averse optimal control}
\subsection{Measuring risk} \label{sec:measuring_risk}
Let \(\NN=\N_{[1,n]}\) be a discrete sample space. A probability measure thereon
can be identified by a probability vector \(p\in\Re^{n}\) with \(\sum_{i=1}^n p_i = 1\), \(p_i\geq 0\)
for \(i\in\NN\).
Let \(Z:\NN\to\Re\) be a real-valued random variable on \(\NN\) which represents a random cost; 
for \(i\in\NN\) let \(Z_i = Z(i)\). The vector \((Z_{i})_{i{\in}\NN}\) identifies
the random variable \(Z\).

The \textit{expectation} of a random variable \(Z\) with respect to the 
probability vector \(p\) is defined as 
\begin{align}
 \E_p[Z] \equiv \E_p[Z(i); i] = \sum_{i\in\NN}p_{i}Z_{i}.
\end{align}
The notation \(\E_p[Z; i]\) is to emphasize that the expectation is taken
with respect to \(i\).

A \textit{risk measure} on \(\Re^n\) is a mapping \(\rho:\Re^{n}\to\Re\). 
It is called \textit{coherent} if it satisfies 
the following properties~\cite[Sec.~6.3]{shapirolectures} for 
\(Z,Z'\in\Re^n\), \(\alpha\geq 0\), \(\lambda\in[0,1]\)
\begin{enumerate}
 \item[A1.] \textit{Convexity}. 
	    \(\rho(\lambda Z + (1-\lambda) Z') 
	      \leq \lambda \rho(Z) + (1-\lambda) \rho(Z')\),
 \item[A2.] \textit{Monotonicity}. 
	      \(\rho(Z)\leq \rho(Z')\) whenever \(Z\leq Z'\),
 \item[A3.] \textit{Translation equivariance}. 
	    \(\rho(c1_n+Z)=c+\rho(Z)\),
 \item[A4.] \textit{Positive homogeneity}. 	     
	      \(\rho(\alpha Z) = \alpha \rho(Z)\).
\end{enumerate}
Trivially, the \textit{expectation} is a coherent risk measure 
and so is the \textit{essential maximum} \(\essmax[Z] {}\dfn{} \max\set{Z_i}[p_i>0, i\in \NN]\). 
A popular risk measure is the \textit{average value-at-risk}, also known as 
\textit{conditional value-at-risk} or \textit{expected shortfall}, which is defined as 
\begin{equation*} 
	\AVAR_\alpha[Z] 
{}={}
	\begin{cases}
	    \min\limits_{t\in\Re}\, \set{t+\alpha^{-1}\E_p\plus{Z-t}},
		&\hspace*{-5pt}\alpha{}\in{}(0,1]
	\\
	    \essmax(Z),
		&\hspace*{-5pt}\alpha=0.
	\end{cases}
\end{equation*}

As a result of assumptions A1--A4, coherent risk measures can be written in 
the following dual form~\cite[Thm.~6.5]{shapirolectures}
%
%
\begin{equation}%
\label{eq:rho_duality_expectation}
	\rho[Z]
{}={} 
	\max_{\mu \in \mathcal{A}(p)}  \E_\mu[Z],
\end{equation}
where \(\mathcal{A}(p)\subseteq\Re^n\) is a compact convex set 
of probability vectors containing \(p\) which we shall call 
the \textit{ambiguity set} of \(\rho\).
We may think of a coherent risk measure as the worst-case expectation
with respect to a probability distribution taken from a  set of 
probability vectors.
We call \(\rho\) a \textit{polytopic} risk measure if \(\mathcal{A}(p)\) is a
polytope, i.e., it can be described by 
\(\rho(Z) = \max\set{\mu^{\tttop} Z}[1_{n}^{\tttop}\mu=1, F(p)\mu \leq b(p)]\)
for some \(F(p)\in\Re^{q\times n}\) and \(b(p)\in\Re^{q}\).
The expectation, the essential maximum and \(\AVARa\) are polytopic risk measures.
The ambiguity set of \(\AVAR_{\alpha}\) for \(\alpha\in [0,1]\) is the polytope
\(
 \mathcal{A}_\alpha(p) {}={} \set{
	\mu\in\Re^n }[
	    \textstyle\sum_{i=1}^n \mu_i =1,\, \mu_i \geq 0, \alpha \mu_i \leq p_i{}].
\)
The ambiguity set \(\mathcal{A}_0(p)\) is the whole probability simplex.
Apparently \(\AVARa\) is a polytopic risk measure.
\(\AVARa\) interpolates between  the risk-neutral expectation operator (\(\AVAR_1=\E_p\), 
with \(\mathcal{A}_0(p)=\set{p}\)) and the worst-case essential maximum (\(\AVAR_0=\essmax\)).%
%
%
%
\subsection{Markovian switching systems} \label{sec:markovian_switching_systems}
In this work we consider Markovian switching systems
\begin{equation}
    x_{k+1} 
{}={} 
    f(x_k, u_k, i_k)
 \label{eq:markovian-system},%
\end{equation}
driven by the random parameter \(i_k\) which is a time-homogeneous
Markov chain with values in a finite set \(\NN=\N_{[1,n]}\)
with transition matrix \(P=(p_{ij})\in\Re^{n\times n}\), that is
\(\probability[i_{k+1}=j\mid i_k=i]=p_{ij}\)~\cite{mjls2005costa}.
We call the states of this Markov chain, the \textit{modes} of~\eqref{eq:markovian-system}.
We denote the \textit{cover} of each mode by \(\cov(i)\dfn \set{j\in\NN}[p_{ij}>0]\).
We assume that at time \(k\) we measure the full state \(x_k\) and the value of \(i_k\).
As the probabilistic information available up to time \(k\) is fully described 
by the pair \((x_k, i_k)\), the control actions \(u_k\) may be decided by 
a causal control law  \(u_k = \kappa_k(x_k, i_k)\). 
This formulation aligns with that of the classic textbook~\cite{mjls2005costa},
but there exist formulations where 
\(i_k\) is not known at time \(k\) and the control law is a function
of \(x_k\) only~\cite{chow2017arxiv}.

Each 
\(
	f({}\cdot{},{}\cdot{},i)
{}:{}
	\Re^{n_x}\times \Re^{n_u}
{}\to{} 
	\Re^{n_x}
\), 
\(
	i
{}\in{}
	\NN
\), 
is assumed to be continuous and satisfy \(f(0,0,i)=0\).
\ac{mjls} are a special case of~\eqref{eq:markovian-system} with 
\(f(x,u,i) = A_{i} x + B_{i} u\), \(i\in\NN\).
System~\eqref{eq:markovian-system} is subject to the
joint state-input  constraints
\begin{equation}\label{eq:constraints}
 (x_k, u_k) \in Y_{i_k},
\end{equation}
and we shall assume that for all \(i\in\NN\), \(Y_{i}\)
are nonempty, closed sets containing the origin.%
%
\subsection{Markov risk measures}\label{sec:crm}
Consider the space of pairs \((i,j)\) in
\(
    \Omega
{}\dfn{}
    \NN\times\NN
\) 
equipped with the \(\sigma\)-algebra \(\mathscr{F}=2^\Omega\)
and the probability measure 
\(
    \probability[\set{(i,j)}]
{}={}
    p_{ij}
\).
The conditional probability conditioned by the knowledge of \(i\)
can be identified with the probability vector \(P_{i}\) 
--- the \(i\)-th row of \(P\).
For a random variable \(Z:\Omega\to\Re\), the conditional expectation of \(Z\) conditioned by 
\(i\), denoted as \(\E_{i}[Z;j]\), is a random variable on \(\NN\), that is 
\(
    \NN
{}\ni{} 
    i
{}\mapsto{}
    \E_{i}[Z;j]
{}\in{}
    \Re
\), with 
\begin{align}\label{eq:conditional-expectation}
 \E_{i}[Z;j]
{}\dfn{}
  \E_{P_{i}}[Z;j]
{}={}
  \hspace{-0.2em}
  \textstyle\sum\limits_{j{}\in{}\NN}
    p_{ij}Z(i,j).
\end{align}
%
%
We may extend this definition to define conditional 
variants of risk measures. Following~\eqref{eq:conditional-expectation},
we give the following definition

\begin{definition}[Markov risk measure]
Given a coherent risk measure \(\rho\) with ambiguity set \(\mathcal{A}\) and a 
probability transition matrix \(P\) of a Markov chain, we define the 
Markov risk measure \(\Crho[i]{Z}[][j]\) as
\begin{align}
   \Crho[i]{Z}[][j]
{}={}
      \max_{\mu\in\mathcal{A}(P_{i})}
      \smashunderbracket{
	\textstyle\sum_{j\in\NN}
	    \mu_{j} Z(i,j)
      }
      {\E_{\mu}[Z;j]},
\end{align} 
for all random variables \(Z:\Omega\to\Re\).
\end{definition}

This definition falls into the general framework of~\cite{Ruszcz2010}.
This way, with every \(i\) we associate the coherent
risk measure \(\Crho[i]{Z}[][j]\).
As with the expectation, the notation \(\Crho[i]{Z}[][j]\) is to emphasize
that the risk is computed with respect to \(j\).%
%
\subsection{Risk-averse optimal control and \acl{dp}}
\label{sec:risk-averse-oc-and-dp}

Consider a \textit{stage cost} function \(\ell\in\fcns(\Re^{n_x}\times\Re^{n_u},\NN)\)
and a \textit{terminal cost} \(\ell_N \in \fcns(\Re^{n_x},\NN)\).
Functions \(\ell\) are extended-real-valued, therefore, they can encode constraints 
such as~\eqref{eq:constraints} by taking \(\dom \ell(\cdot,\cdot,i) = Y_{i}\), \(i\in\NN\).
Likewise, \(\ell_N\) can encode constraints on the terminal state of the form 
\(x_N\in X^f_{i_N}\) by taking \(\dom \ell_{N}(\cdot,i) = X^{f}_{i}\), \(i\in\NN\),
where \(X^f_{i}\) contain the origin in their interiors. 
We may now introduce the following finite-horizon risk-averse optimal control problem 
\begin{align}\label{eq:risk-averse-problem}
	\minimize_{
	    \substack{ 
		u_0
	    } 
	}
{\,}&{}    
	\ell(x_0, u_0,i_0)  
{}+{}
	\rho_{i_0}
	\bigg[{}\inf_{u_1} \ell(x_1, u_1,i_1)  
\notag
\\
&{}+{}
	  \rho_{i_1}
	  \Big[
	      {}\inf_{u_2}{}\ell(x_2,u_2,i_2)  + \cdots 
\notag
\\        
{}&{}+ \rho_{i_{N{-}1}}[\ell_{N}(x_N{,}i_N){;}i_N]\cdots
	  ;i_2\Big]
	;i_1\bigg],
\end{align}
where \(x_{k+1} = f(x_k, u_k, i_k)\), for all \(k\in\N_{[0,N-1]}\).
As it will become evident in what follows, each one of the infima at stage \(k\) in
\eqref{eq:risk-averse-problem} is parametric in \(x_k\) and \(i_k\), that is,
the minimization takes place over causal control laws \(u_0,\ldots,u_{N-1}\).
Note that under assumptions A1 and A2, we may interchange the Markov risk measures 
with the infima~\cite[Prop.~6.60]{shapirolectures} leading 
to risk-averse multistage formulations discussed in~\cite[Sec.~6.8.4]{shapirolectures}.

Problem~\eqref{eq:risk-averse-problem} can be described by a \ac{dp} 
recursion.
Inspired by~\cite[Sec.~6.8]{shapirolectures}, for a \(V \in \fcns(\Re^{n_x}, \NN)\) 
we define the \ac{dp} operator
\(\T:\fcns(\Re^{n_x}, \NN) \to \fcns(\Re^{n_x}, \NN)\) so that 
\begin{align*}
 (\T V)(x, i) 
    &{=} \inf_{u} \set{\ell(x, u, i) 
	  {+}\rho_{i}\left[V(f(x,u,i),j);j\right]}\\
	  &{=}\inf_{u} \ell(x, u, i)
+ 
	\max_{\mu\in\mathcal{A}(P_i)} 
	  \sum_{j\in\NN} \mu_{j} V(f(x,u,i),j).
\end{align*}%
Let \((\S V)(x, i)\) be the corresponding set of minimizers for the 
optimization problem involved in the definition of \((\T V)(x,i)\). 
This defines the following DP recursion%
\begin{subequations}\label{eq:dp-recursion}%
 \begin{align}
  V_{k+1}^\star &= \T V_k^\star,\\
  \U_{k+1}^\star &= \S V_k^\star,
 \end{align}
\end{subequations}
for \(k\in\N_{[0,N-1]}\) with \(V_0^\star(x,i) \dfn \ell_{N}(x,i)\), \(i\in\NN\).
For 
\(
	C
{}={}
	\seq{C_{i}}[i\in\NN]
\)
with 
\(
	C_{i}
\subseteq
	\Re^{n_x},
\)
we define the mode-dependent \textit{predecessor operator} 
\(
	R(C)
{}={}
	\seq{R_i(C)}[i\in\NN]
\) 
with 
\(
    R_{i}(C)
{}={}
  \set
    {x{}\in{}\Re^{n_x}}
    [
      \exists u {}\in{} \Re^{n_u}, 
      (x,u) {}\in{} Y_i, 
      f(x,u,i) {}\in{} \bigcap_{j\in\cov(i)} C_{j}
    ]
\).
Next, we present some fundamental properties of the DP operator~\(\T\).

\begin{proposition}\label{prop:dp_properties}
If \(\ell_{N}({}\cdot{}, i)\) are proper, \ac{lsc} and \(\ell(\cdot,\cdot,i)\) are 
proper, \ac{lsc} and level bounded in \(u\) locally uniformly in \(x\) for all \(i\in\NN\),
then for all \(i\in\NN\):
\begin{enumerate*}[label=(\roman*)]%
 \item \(\T V \in \fcns(\Re^{n_x},\NN)\) for \(V\in\fcns(\Re^{n_x},\NN)\),             
 \item \(V_k^\star({}\cdot{},i)\) are \ac{lsc},                                        
 \item \(\dom V_k^\star({}\cdot{},i)\) \(= \dom \U_k^\star({}\cdot{},i) \neq \emptyset\),  
 \item \(\U_k^\star\) is compact-valued,                                               
 \item \(\dom(V_{k+1}^{\star}){=}R(\dom(V_{k}^{\star}))\).                             
\end{enumerate*}%
\begin{pf}%
The proof goes along the lines of~\cite[Thm.~11a]{smpc4mss} 
using~\cite[Prop.~1.17, Prop.~1.26(a)]{rockafellar2011variational}.
\end{pf}%
\end{proposition}%
We may easily verify the monotonicity property
\(
  \T V \leq \T V' ,\ \textrm{for all }\ V,V'\ \textrm{with }\ V\leq V',
\)
following~\cite{bertsekas-book-vol2}. 
An observation that will prove useful in what follows is that if 
\(\T\ell_N \leq \ell_N\), then \(V_{k+1}^\star \leq V_k^\star\).
The above risk-averse optimal control problem leads naturally to 
the statement of a risk-averse \ac{mpc} problem
where control actions are computed by a control law 
\(\kappa_{N}^\star(x,i)\in \U_{N}^\star(x,i)\).
In Section~\ref{sec:risk-averse-stability} we state an appropriate
risk-based notion of stability and provide conditions on \(\ell_N\) for 
the \ac{mpc}-controlled system \(x_{k+1}=f(x_k, \kappa_N^\star(x_k,i_k),i_k)\)
to be stable.%

\section{Risk-averse stability}%
\label{sec:risk-averse-stability}
Consider the following Markovian switching system which is controlled
by some control law \(u_k = \kappa(x_k,i_k)\)%
\begin{align}\label{eq:system-unactuated}%
 x_{k+1} = f^{\kappa}(x_k,i_k) \dfn f(x_k,\kappa(x_k,i_k), i_k),
\end{align}%
subject to the constraints 
\((x_k, i_k) \in X^\kappa \dfn \set{(x,i)}[(x,\kappa(x,i))\in Y_{i}]\).
For convenience, we introduce the notation \(X^\kappa_{i} = \set{x}[(x,\kappa(x,i))\in Y_{i}]\),
for \(i\in\NN\).
Let 
\(
	i_{[k]}
{}={}
	(i_0,i_1,\ldots,i_{k})
\) 
denote an  admissible path of length \(k\) of the Markov chain \(\seq{i_t}[t\in\N]\),
that is, \(i_{t+1}\in\cov(i_t)\) for  \(t\in\N_{[0,k-1]}\).
For a given initial state \(x_0\), the solution of~\eqref{eq:system-unactuated} 
at time \(k\) is denoted as \(\phi(k,x_0,i_{[k-1]})\).

In order to be able to define risk-based notions of stability, 
we must first introduce an appropriate notion of invariance for 
Markovian switching systems~\cite{smpc4mss}.

\begin{definition}[Uniform invariance]
 Let \(X=\seq{X_{i}}[i\in\NN]\) be a 
 collection of nonempty closed subsets of \(\Re^{n_x}\) 
 and \(X_{i} \subseteq X^\kappa_{i}\). 
 \(X\) is called \ac{ui} for~\eqref{eq:system-unactuated}
 subject to  constraints 
\(
	x 
{}\in{}
	X^\kappa_{i}
\)	
if 
\(
	f^{\kappa}(x,i) 
{}\in{} 
	\bigcap_{j\in \cov(i)}
	X_{j},
\)
whenever \(x {}\in{}  X_{i}\) for all \(i\in\NN\).
\end{definition}

For the controlled system~\eqref{eq:system-unactuated}, the predecessor operator 
is now defined as 
\(
    R_{i}(C)
{}={}
  \set
    {x{}\in{}X^\kappa}
    [f^{\kappa}(x,i){}\in{}\bigcap_{j\in\cov(i)}C_{j}]
\).
We have that \(C\) is \ac{ui} if and only if 
\(C_{i} \subseteq R_{i}(C)\) for all \(i\in\NN\)~\cite{smpc4mss}.

Given a coherent risk measure \(\rho\) and a random variable 
\(\psi(i_0, i_1, \ldots, i_k)\), let  \(\bar{\rho}_1[\psi] = 
\Crho[i_0]{\psi(i_0, i_1,\ldots, i_k)}[][i_1]\)
and recursively define 
\(
      \bar{\rho}_{k} 
{}={}
      \bar{\rho}_{k-1} 
{}\circ{}
      \Crho
	[i_{k-1}]
	{\cdot}
	[]
	[i_{k}]
\),
that is 
\(
	\bar{\rho}_{k}[\psi] 
= 
	\Crho
	    [i_0]
	    {
	    \Crho
	      [i_1]
	      {
	      \cdots\Crho[i_{k-1}]{\psi(i_0, i_1,\ldots, i_k)}[][i_k]\cdots
	      }
	      [][i_2]
	    }
	    [][i_1]
\)
\cite[Sec.~6.8.2]{shapirolectures}.

We may now give the following stability notion~\cite{chow-risk-averseMPC}.

\begin{definition}[Risk-square exponential stability]\label{def:rss}
We say that the origin is \ac{rses} for system~\eqref{eq:system-unactuated}
over a set 
\(
	X
{}={} 
	\seq{X_i}[i\in\NN]
  \) 
if \(X\) is \ac{ui} and for \(x_0{}\in{}X_{i_0}\)
\begin{equation*}
    \bar\rho_{k-1} \bigl[
	  \|\phi(k,x_0,i_{[k-1]})\|^2
	  \bigr]
{}\leq{}
    \lambda \beta^{k} \|x_0\|^2,
\end{equation*}
for all \(k\in\N\), for some \(\beta{}\in{}[0,1)\), \(\lambda {}\ge{} 0\).
\end{definition}
\ac{rses} entails that the origin is exponentially mean-square stable for system%
~\eqref{eq:system-unactuated} not only for the nominal probability distribution, 
but also for those probability distributions in the ambiguity set of the risk
measure. 
In the unconstrained case, \ac{rses} corresponds to the notion of uniform 
global risk-sensitive exponential stability which is defined using the 
notion of dynamic risk measures~\cite{chow-risk-averseMPC}.
If the underlying risk measure is the expectation operator, then 
\ac{rses} reduces to mean-square exponential stability, whereas, if it is the 
essential supremum operator, it yields the definition of robust exponential 
stability. Additionally, since all coherent risk measures are lower bounded 
by the expectation, \ac{rses} is a stronger notion of stability compared to mean-square
stability.
The following lemma provides Lyapunov-type stability conditions
for \ac{rses}.

%
%
\begin{lem}[\ac{rses} conditions]\label{lem:mrses}
 Suppose there is a \(V\in\fcns(\Re^{n_x},\NN)\), proper, \ac{lsc} function such that
 \begin{enumerate}[label=(\roman*)]
  \item \label{item:domV_assumptions}\(\dom V\) is a \ac{ui} set
  \item \label{item:lyap_condition}
	\(
	    \Crho[i]{V(f^\kappa(x,i),j)}[][j]
	{}-{}
	     V(x,i) 
	{}\leq{}
	     -c \|x\|^2
	\), 
        for some \(c > 0\) 
        for all \((x,i)\in \dom V\).
 \end{enumerate}%
Then, 
\(
    \bar{\rho}_{k}
      \bigl[
	  \textstyle\sum_{t=0}^{k-1} \|\phi(t,x_0,i_{[t-1]})\|^2
      \bigr],
\) 
is uniformly bounded in \(k\) for \((x_0, i_0)\in \dom V\).
If, additionally, 
\begin{enumerate}[resume,label=(\roman*)]
 \item \label{item:quad_condition}
    for all \((x,i) {\in} \dom V\), 
    \(
	\alpha_1 \|x\|^2 
{}\leq{}
	V(x,i) 
{}\leq{}
	\alpha_2 \|x\|^2
    \), 
for some \(\alpha_1, \alpha_2>0\),
\end{enumerate}
then, the origin is \ac{rses} for system~\eqref{eq:system-unactuated} over \(\dom V\).
%
%
\begin{pf}
The proof can be found in the appendix.
\end{pf}
\end{lem} 
The uniform boundedness condition in Lemma~\ref{lem:mrses} is reminiscent of the notion
of stochastic stability in~\cite[Sec.~3.3.1]{mjls2005costa}.
In fact, if the risk measure in Lemma~\ref{lem:mrses} is the expectation operator,
then the uniform boundedness condition is equivalent to mean-square 
stability~\cite[Thm.~3.9(6)]{mjls2005costa}.

We call a function 
\(
	V
{}\in{}
	\fcns(\Re^{n_{x}},\NN)
\) 
which satisfies all requirements of 
Lemma~\ref{lem:mrses}, a (mode-dependent) \textit{risk-averse Lyapunov function}.
We may now state conditions on the stage cost \(\ell\) and the 
terminal cost \(\ell_N\) which
entail \ac{rses} for the risk-averse \ac{mpc}-controlled system.%

\section{Risk-averse MPC}
\subsection{Risk-averse MPC stability}
%
%
\begin{thm}[\ac{rses} of \ac{mpc}]\label{thm:mrses-cond-tv}
Suppose that 
\begin{enumerate*}[label=(\roman*)]
 \item \label{it:mpc-mrses-ell-bound}  \(c\|x\|^2\) \(\leq  \ell(x,u,i)\) for some \(c>0\) for all \((x,u) \in Y_{i}\), \(i\in\NN\)
 \item \label{it:mpc-mrses-ellN-bound} \(\ell_{N}(x,i) \leq d \|x\|^2,\) for some \(d  > 0\) for all \(x\in X^{f}_{i}\),
 \item \label{it:mpc-mrses-Xfi-origin} \(X^f_i\) contain the origin in their interiors 
 \item \label{it:mpc-mrses-VNstar-loc-bounded} \(V_N^\star\) is locally bounded over its domain, that is, for every compact set \(\bar{X}\subseteq \dom V_N^\star\),
      there is an \(M{}\geq{}0\) so that \(V_N^\star(x,i)\leq M\) for all 
     \((x,i)\in \bar{X}\) and%
\end{enumerate*}
\begin{equation}\label{eq:tv-leq-v}%
\T \ell_N \le \ell_N.
\end{equation}%
Then, the origin is \ac{rses} for the risk-averse 
\ac{mpc}-controlled system
\(
	x_{k+1}
{}={}
	f(x_k, \kappa_N^\star(x_k, i_k),i_k)
\)
over all compact uniformly invariant subsets of \(\dom V_N^\star\).
%
%
\begin{pf}%
The proof can be found in the appendix.
\end{pf}
\end{thm} 

In Thm.~\ref{thm:mrses-cond-tv} we show that \(V_N^\star\) is a mode-dependent risk-averse
Lyapunov function over compact uniformly invariant subsets of \(\dom V_N^\star\). 
We shall use this result in the following sections 
to design risk-averse stabilizing \ac{mpc} controllers for \ac{mjls} as well as 
nonlinear Markovian switching systems.
Note that Condition \ref{it:mpc-mrses-VNstar-loc-bounded} in Thm. \ref{thm:mrses-cond-tv} holds if
the following assumption is satisfied (see \cite[Prop. 2.15]{rawlings2009model})
\begin{assumption}[Local boundedness of \(V_N^\star\)]\label{as:VN-loc-bounded}
 For~all \(i\in\NN\), functions \(\ell(\cdot, \cdot, i)\) and \(\ell_N(\cdot, i)\) 
 are continuous on their domains, and
  the sets \(U_{i}(x)\dfn\{u\in\Re^{n_u} {}\mid{} (x,u) \in Y_i\}\)
 are compact and bounded uniformly in \(x\).
\end{assumption}
Additionally, because of the monotonicity property of \(\T\) and since \(\T \ell_N \leq \ell_N\), 
condition~\eqref{eq:tv-leq-v} implies  \(V_{k+1}^\star\leq V_k^\star\),
thus 
\(
	\dom(V_{k}^\star) 
{}\subseteq{} 
	\dom(V_{k+1}^\star)
{}={}
	R(\dom V_k^\star)
\) 
(Prop.~\ref{prop:dp_properties}), thus \(\dom V_k^\star\) is \ac{ui}.

\subsection{Risk-averse \ac{mpc} design for \ac{mjls}}
\label{sec:risk-averse-mjls}
Here we provide \ac{rses} conditions and design
guidelines for risk-averse \ac{mpc} of \ac{mjls}~\cite{mjls2005costa}, 
that is \(f(x,u,i) = A_{i} x + B_{i} u\), using a quadratic stage cost 
\(
      \ell(x,u,i) 
{}={}
      x^{\tttop}Q_{i} x 
{}+{}
      u^{\tttop}R_{i} u 
{}+{} 
      \delta_{Y_{i}}(x,u), 
\)
with \(Q_{i} \in \Sp^{n_x}\), \(R_{i} \in \Spp^{n_u}\) and 
\(Y_{i}\) are polytopes with the 
origin in their interiors. The terminal cost function is taken to be 
\(
      \ell_{N}(x,i) 
{}={} 
      x^{\tttop}P^{f}_{i} x 
{}+{} 
      \delta{\scriptscriptstyle{ X^f_{i}}}(x)
\)
with \(P^{f}_{i} {}\in{} \Spp^{n_x}\) and \(X^{f}_{i}\).
We shall derive conditions on \(P^{f}_{i}\) and \(X^{f}_{i}\) so that the stabilizing 
conditions of Thm.~\ref{thm:mrses-cond-tv} are satisfied.
Condition \(\T \ell_N \leq \ell_N\) is equivalent to%
\begin{subequations}%
\begin{align}%
&\min_{u}
\set{
 x^{\tttop}Q_{i}x 
    + u^{\tttop}R_{i}u\notag \\    
&
  \quad+\Crho
      [i]
      {x^{\tplus\tttop}P^f_{j}x^{\tplus}}
      []
      [j]    
 } 
	  {}\le{} 
      x^\tttop P^{f}_{i}x,\label{eq:xyz-1}\\
&\dom(\T \ell_N)\supseteq  \dom \ell_N \Leftrightarrow R(X^f)\supseteq X^f,\label{eq:xyz-2}
\end{align}%
\end{subequations}%
where
\(
      x^{\tplus}
{}={}
      f(x,u,i)
\)
and the minimization in~\eqref{eq:xyz-1} is over the space of 
admissible causal control laws \(u = \kappa(x, i)\) so that 
\((x, i)\in X^\kappa\).
An upper bound to the left hand side of~\eqref{eq:xyz-1} 
is obtained by parametrizing
\(
u = K_{i} x.
\)
We introduce the shorthand notation  
\(\bar{A}_{i} = A_{i} + B_{i}K_{i}\) 
and
\(\bar{Q}_{i} = Q_{i} + K_{i}^{\tttop}R_{i}K_{i}\),
for \({i}\in\NN\).
%
Condition~\eqref{eq:xyz-2} means that \(X^f\) is a \ac{ui} set for the system 
\(
        x_{k+1} 
{}={}
        (A_{i_k}+B_{i_k}K_{i_k})x_k
\)
under the prescribed constraints.
Such a set can be determined by the fixed-point iteration
\(\mathcal{O}_{k+1} = R(\mathcal{O}_k)\) with \(\mathcal{O}_0 =\set{(x,i)}[(x,K_i x)\in Y_i]\).
If this iteration converges in a finite number of iterations --- a sufficient
condition for which is given in~\cite[Lem. 21]{smpc4mss} --- to a set 
\(\mathcal{O}_{\infty}\), this is a \textit{polytopic} \ac{ui} set.

Assuming that \(\rho\) is a polytopic Markov risk
measure with ambiguity set 
\(
	  \mathcal{A}(P_{i}) 
{}={}
	  \conv\set{\mu_{i}^{\scriptscriptstyle(l)}}_{l\in\N_{[1,s_{i}]}}
\) and 
using its dual representation, condition~\eqref{eq:xyz-1} becomes
\(
	\bar{Q}_{i} 
+
	\sum_{j\in \cov(i)} 
	\mu_{ij}^{\scriptscriptstyle(l)}
\)
\(
	(\bar{A}_{i}^{\tttop} P^f_{j} \bar{A}_{i} ) 
\preccurlyeq  
	P^f_{i}
\)
for all \(i\in\NN\) and \(l\in\N_{[1,s_{i}]}\).
This condition can be cast 
as a \ac{lmi} by a change of variables
\(
	(P^f_{i})^{-1} 
=
	M_{i}
\), 
\(
	K_{i} 
=
	Y_{i} M_{i}^{-1}
\), 
\(
	F_{i}^{l} 
= 
	\smallmat{ \sqrt{\mu_{i{}1}^{\scriptscriptstyle(l)}}I \ldots  \sqrt{\mu_{in}^{\scriptscriptstyle(l)}I} }
\)
and \(M = \blkdiag(M_1,\ldots, M_{n})\):%
%
\begin{equation}\label{eq:LMI-implication}%
 \smallmat{
  M_{i}   & (A_{i} M_{i} + B_{i} Y_{i})^{\tttop} F_{i}^{l} 
                 & M_{i} Q_{i}^{1/2} 
                          & Y_{i}^{\tttop} R_{i}^{1/2}\\
 \ast & M    & 0       & 0\\
 \ast & \ast & I       & 0\\
 \ast & \ast & \ast    & I 
 }
 \succcurlyeq 0,
\end{equation}%
for all \(i\in\NN\) and \(l\in\N_{[1,s_{i}]}\).
The left hand side of \eqref{eq:LMI-implication} is a symmetric matrix, therefore,
we show only its upper block triangular part and replaced the lower block triangular 
part by asterisks (\(*\)) to simplify the notation.
Solving this \ac{lmi} for \(M_{i}\in\mathcal{S}_{\tplus}^{n_x}\) and \(Y_{i}\in\Re^{n_u\times n_x}\)
yields the linear gains \(K_{i}\) and the cost matrices \(P^f_{i}\).
\ac{lmi}~\eqref{eq:LMI-implication} has to be solved once offline to determine matrices \(P^f_i\).%

\subsection{Risk-averse \ac{mpc} design for nonlinear Markovian switching systems}
\label{sec:linearization-design}

For nonlinear systems, an obvious choice for the terminal cost function
would be \(\ell_{N}(x,i) = \delta_{\set{0}}(x)\) --- meaning, 
\(X^{f}_{i} = \set{0}\) for \(i\in\NN\) --- 
but that would lead to a very conservative design.
Here we exploit the system linearization at the origin
to determine a terminal cost function and terminal constraints which 
render the \ac{mpc}-controlled system \ac{rses}.
We shall first draw the following assumption for the nonlinear dynamics:
\begin{assumption}\label{as:f-beta-smooth}
 For each \(i\in\NN\), \(f({}\cdot{},{}\cdot{},i)\) is differentiable 
 with \(L_{i}\)-Lipschitz Jacobian.
\end{assumption}%
We use a parametric controller
of the form \(\kappa(x,i) = K_{i} x\) and define the associated 
closed-loop function \(f^\kappa(x,i) = f(x, K_{i} x, i)\), \(i\in\NN\).
Function \(f^\kappa({}\cdot{},{}\cdot{},i)\) can be written as a 
composition of \(f({}\cdot{},{}\cdot{},i)\) with the linear
mapping \(W_{i}:(x,u)\mapsto (x, K_i x)\), therefore, its 
Jacobian matrix will be Lipschitz-continuous with Lipschitz constant
\(L_{i} \|W_{i}\|^2\) which is bounded above by%
\begin{equation}%
	  \beta_{i} 
{}\dfn{} 
	  L_{i} (1 {}+{} \|K_{i}\|^2).
\end{equation}%
The linearization of the nonlinear system at the origin 
is an \ac{mjls} \(x_{k+1} = \hat{f}(x_k,u_k,i_k) \dfn A_{i_k} x_k + B_{i_k} u_k\)
with \(A_{i_k}\) and \(B_{i_k}\) given by the Jacobian matrices,
with respect to \(x\) and \(u\) respectively, of \(f\) at the origin.
That is, \(A_{i} = J_xf(0,0,i)\),  \(B_{i} = J_u f(0,0,i)\).
For notational convenience, we define the following quantities%
\begin{align*}
	\hat{f}^\kappa(x,i) 
{}\dfn{}& 
	(A_i + B_i K_i)x,
\\
	\mathcal{L}\ell_N(x,i)  
{}\dfn{}&
	\Crho[i]{\ell_N(f^{\kappa}(x,i),j);j} {}-{} \ell_{N}(x,i),
\\
	\mathcal{L}'\ell_N(x,i) 
{}\dfn{}& 
	\Crho[i]{\ell_N(\hat{f}^{\kappa}(x,i),j);j} {}-{} \ell_{N}(x,i),
\\
	\Delta(x,i) 
{}\dfn{}& 
	\mathcal{L}\ell_N(x,i) {}-{} \mathcal{L}'\ell_N(x,i).
\end{align*}%
The objective is to design the terminal cost and terminal 
constraints for the risk-averse \ac{mpc} problem using \(\mathcal{L}'\ell_N\)
to yield an \ac{lmi}. 
While our design will be based on the linearized dynamics,
we need to account for the linearization error.
To this end, we shall derive a quadratic upper bound for 
\(|\Delta(x_k,i_k)|\) in a neighborhood of 
the origin.

\begin{thm}\label{thm:linearization}
 Suppose that Assumptions~\ref{as:VN-loc-bounded} and~\ref{as:f-beta-smooth} hold and%
\begin{equation}\label{eq:linearization-1}%
  \mathcal{L}'\ell_N(x,i) \leq - x^{\tttop}(\bar{Q}_{i} + m_{i}I)x,
\end{equation}%
for \(i\in\NN\), \(m_{i}>0\), 
\(\ell\), \(\ell_N\) and \(X^f\) satisfy the requirements of Thm. \ref{thm:mrses-cond-tv}
with \(X^f_i \subseteq \mathcal{B}_{\delta_i}\) for some \(\delta_i>0\) for all \(i\in\NN\), 
\begin{equation}\label{eq:delta}%
 \sigma_i \dfn \max_{j\in\cov(i)}\|P^f_{j}\|
 \,
  (\tfrac{\beta_{i}^2 \delta_{i}^2}{4} +\beta_{i}\|\bar{A}_{i}\| \delta_{i}) <  m_{i},
\end{equation}%
and \(\ell(x, K_i x, i) \leq x^{\tttop} (\bar{Q}_i + (m_i-\sigma_i)I)x\).
If \(X^f\) is a \ac{ui} set for \eqref{eq:system-unactuated}, then the origin is 
\ac{rses} for the \ac{mpc}-controlled system 
\(
	x_{k+1} 
{}={}
	f(x_k, \kappa_N^\star(x_k, i_k), i_k)
\)
over the compact \ac{ui} subsets of \(\dom V_N^{\star}\).

\begin{pf}
The proof can be found in the appendix.
\end{pf}
\end{thm}%

According to Thm.~\ref{thm:linearization}, one first needs
to select \(m_{i}>0\) for each \(i\in\NN\) such that~\eqref{eq:linearization-1}
holds true. In the most common case where \(\ell\) and \(\ell_N\) are quadratic 
functions, this is precisely an \ac{lmi} of the form~\eqref{eq:LMI-implication}
with \(Q_{i}+m_{i} I\) in place of \(Q_{i}\) solving which
we obtain matrices \(K_{i}\) and \(P^{f}_{i}\) and determine the constants 
\(\beta_{i}\) and find \(\delta_{i}>0\) so that~\eqref{eq:delta} holds.
The last step is to determine a \ac{ui} set \(X^f\) 
for the nonlinear system 
\(
	x_{k+1} 
{}={} 
	f^{\kappa}(x_k, i_k)
\).
We may cast the nonlinear system as a linear one with bounded additive
disturbance 
\(
	x_{k+1} 
{}={} 
	\bar{A}_{i_k} x_k 
{}+{} 
	e(x_k, i_k)
\) --- indeed, as we show in the proof of Thm.~\ref{thm:linearization},
\(
	\|e(x,i)\| 
{}\leq{} 
	\nicefrac{\beta_{i}}{2}\|x\|^2
\).
We may follow the approach of~\cite{SchaichCannon15} in order to 
determine a polytopic robustly invariant set.%

\section{Computationally tractable formulation of risk-averse optimal control problems}
\label{sec:computational}
\newcommand{\node}{\iota}
\newcommand{\nnode}{\eta}

\begin{figure}[t]
 \centering
 \includegraphics[width=0.95\linewidth]{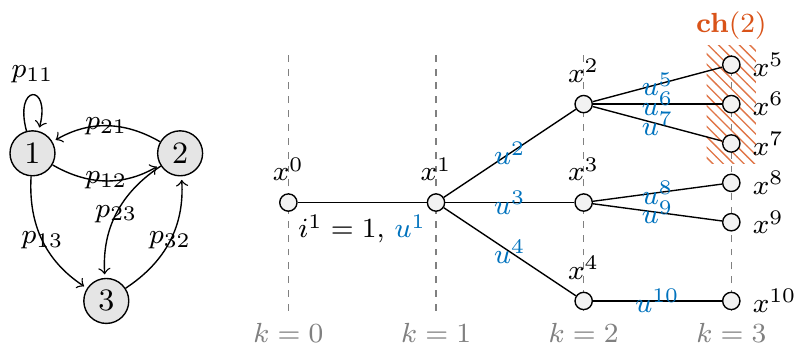}
 \caption{%
      (Left) A Markov chain with three modes and the corresponding transition probabilities,
      (Right) The corresponding tree with \(i_0=1\).
      }
 \label{fig:tree}
\end{figure}

Starting from an initial state \(x_0\) and initial mode 
\(i_0\) and computing control actions according to a causal
control law \(u_k\), the future states of the Markovian system,
up to some future time \(N\), span  a \textit{scenario tree} --- a tree-like structure such as 
the one shown in Fig.~\ref{fig:tree}. Note that the state at a node \(\node\), the 
input and mode leading to that node are denoted as \(x^\node\), \(u^\node\) and \(i^\node\)
respectively.

The possible realizations of the system state at time \(k\)
define the \textit{nodes} of the tree. 
The set of all nodes at stage \(k\) defines the set \(\Omega_k\).
%
The set of nodes in \(\Omega_{k+1}\)
which are reachable from a node \(\node\in\Omega_k\) is called the set of \textit{children} of 
\(\node\) and is denoted by \(\child(\node)\) which is a subset of \(\Omega_{k+1}\). 
The space \(\child(\node)\) becomes a probability space with 
\(
	\probability[\set{\nnode}]
{}={}
	p_{i^{\node}{}i^{\nnode}}
\)	
for \(\nnode {}\in{} \child(\node)\).
As illustrated in Fig.~\ref{fig:tree}, the system dynamics on the scenario tree is 
described by 
\(
	x^{\nnode}
{}={}
	f(x^{\node}, u^{\nnode}, i^{\nnode})
\), 
for 
\(
	\nnode
{}\in{}
	\child(\node)
\)
and \(x^0{}={}x_0\), \(i^1{}={}i_0\).

On the scenario tree, we define a process \(\Phi\) as follows: 
for \(\node{}\in{}\Omega_{N}\) we define 
\(
	\Phi^{\node}
{}\dfn{}
	\Crho
	  [i^{\node}]
	  {\ell_N(x^{\node},i^{\nnode})}
	  []
	  [{\nnode}]
{}={}
	\max_{\mu^{\node}\in\mathcal{A}(P_{i^{\node}})}
	  \sum_{\nnode\in\child(\node)}\mu_{\nnode}^{\node}\ell_N(x^{\node},i^{\nnode})
\).
Moreover, 
\(
	\ell_N(x,i)
\)
\(
=
	\inf_{\ell_N(x,i)\leq \tau}\tau
\). 
When the underlying risk measure is polytopic with 
\(
	\mathcal{A}(p)
=
	\set{\mu\in\Re^n}[\sum_{i=1}^{n}\mu_i=1,F(p)\mu\leq b(p)]
\) with \(b(p)\in\Re^{q}\), then%
\begin{align*}%
	\Phi^{\node}
{}={}&
	\max_{
	    \mu^{\node}\in\mathcal{A}(P_{i^\node})
	    }{} 
	      \inf_{
		  \substack{
			    \ell_N(x^{\node},i^{\nnode})\leq \tau^{\node}_{\nnode},
			    \\ 
			    l {}\in{}\child(\node)
			}
		  }
		  {}
		      \sum_{\nnode\in\child(\node)}
			  \mu_{\nnode}^{\node}\tau^\node_\nnode
\\
{}={}&
	\inf_{
		 \substack{
			    \ell_N(x^{\node},i^{\nnode})\leq \tau^{\node}_{\nnode},
			    \\ 
			    l {}\in{}\child(\node)
			}
	     }{} 
	      \max_{
		      \mu^{\node}\in\mathcal{A}(P_{i^\node})
		    } 
		    {}
			\sum_{\nnode\in\child(\node)}
			    \mu_{\nnode}^{\node}\tau^\node_\nnode
\\
{}={}&
	\inf_{
		\substack
		{
		    \tau^{\node}, y^{\node}\geq 0,{} \lambda^{\node}\in\Re,
		    \\
		    \ell_N(x^{\node},i^{\nnode}) 
		    {}\leq{} 
		    \tau^{\node}_{\nnode},
		    \, 
		    {}
		    l {}\in{}\child(\node)
		    \\
		    \tau^\node {}={} F(P_{i^\node})^{\tttop}y^{\node} + \lambda^{\node} 1_{q}
		}		
	      }
	      {}
		  b(P_{i^\node})^{\ttop}y^{\node} + \lambda^{\node},
\end{align*}%
where in the first equation we interchanged \(\max\) with \(\inf\) using~\cite[Prop.~2.6.4]{BerNedOzd03}
using the fact that the level sets of the mapping 
\(
\tau^\node \mapsto 
\max_{
		      \mu^{\node}\in\mathcal{A}(P_{i^\node})
		    } 
		    {}
			\sum_{\nnode\in\child(\node)}
			    \mu_{\nnode}^{\node}\tau^\node_\nnode
\) are bounded because \(\mathcal{A}(P_{i^\node})\) is compact.
The last equality is because of LP duality.
Traversing indices \(k\) from \(N{-}1\) back to \(1\), we define
\(
\Phi^{\node}
{}\dfn{}
	\Crho
	  [i^{\node}]
	  {\ell(x^{\node}, u^{\nnode}, i^{\nnode}) + \Phi^{\nnode}}
	  []
	  [{\nnode}]
\),
which boils down to%
\begin{align*}
	\Phi^{\node}
{}={} 
\inf_{
		\substack
		{
		    \tau^{\node},{} y^{\node}\geq 0,{} \lambda^{\node}\in\Re
		    \\
		      \ell(x^{\node},u^{\nnode},i^{\nnode}) 
		    {}+{} 
		      \Phi^{\nnode}    
		    {}\leq{} 
		    \tau^{\node}_{\nnode},
		    \, 
		    {}
		    l {}\in{}\child(\node)
		    \\
		    \tau^{\node} {}={} F(P_{i^\node})^{\tttop}y^{\node} + \lambda^{\node} 1_{q}
		}		
	      }
	      b(P_{i^\node})^{\ttop}y^{\node} + \lambda^{\node},
\end{align*}%
for \(\node\in\Omega_k\).
This formulation allows us to deconvolve the 
nested Markov risk measures. Indeed, \(V_N^\star(x_0,i_0)\)
is the optimal value of the following minimization problem%
\begin{align*}%
	\minimize_{
	    \substack{
		x,u,y{}\geq{}0,\lambda,\tau
	    } }{\,}{}
		&{}\ell(x_0,u^1,i_0) 
	{}+{} 
		b(P_{i^1})^{\tttop}y^{1} + \lambda^1\\
	\subjto{\,}{}
		&		{}\ell_N(x^{\node},i^{\nnode})			
			{}\leq{} 
				\tau^{\node}_{\nnode}, \nnode\in\child(\node), \node\in\Omega_N,
		\\	
		&	{}\tau^\node 
		{}={} 
			F(P_{i^\node})^{\tttop}y^{\node} 
		{}+{} 
			\lambda^{\node} 1_{q},
		\\
		& 		{}\ell(x^{\node},u^{\nnode},i^{\nnode}) 
			{}+{} 
				  b(P_{i^\nnode})^{\ttop}y^{\nnode} + \lambda^{\nnode}   
			{}\leq{} 
				\tau^{\node}_{\nnode},			
		\\
		&	{}x^{\nnode} 
			{}={} 
				f(x^{\node}, u^{\nnode}, i^{\nnode}),
		\\
		&{}\nnode\in\child(\node),{}\node\in\Omega_k,{}k\in\N_{[0,N]}.
\end{align*}%
Note that this formulation does not require the enumeration of the vertices 
of \(\mathcal{A}(p)\) which, for instance, in the case of \(\AVARa\) increases
exponentially with the number of modes.
The above optimization problem is solved at every time instant with \(x_0\), \(i_0\) being the 
current state and mode of the system.
Solving this problem yields the optimal control actions \(u^{\node\star}\)
at each node of the scenario tree. The first value, \(u^{1\star}\), defines the risk-averse
MPC controller \(\kappa_N^\star(x,i)=u^{1\star}(x,i)\).
Note that in the particular case of an \ac{mjls} where
stage-wise and terminal costs are quadratic and the constraints
are polyhedral and/or ellipsoidal,
we obtain a \ac{qcqp} which can be solved very efficiently online as we 
show in Section~\ref{sec:examples}.
The above reformulation can be applied to risk measures whose ambiguity 
set is described by a set of conic inequalities (using conic duality)
such as the entropic value-at-risk~\cite{Ahmadi-Javid2012}.%

\section{Illustrative example}
\label{sec:examples}
\begin{figure}[ht]
 \centering
 \includegraphics[width=0.46\columnwidth]{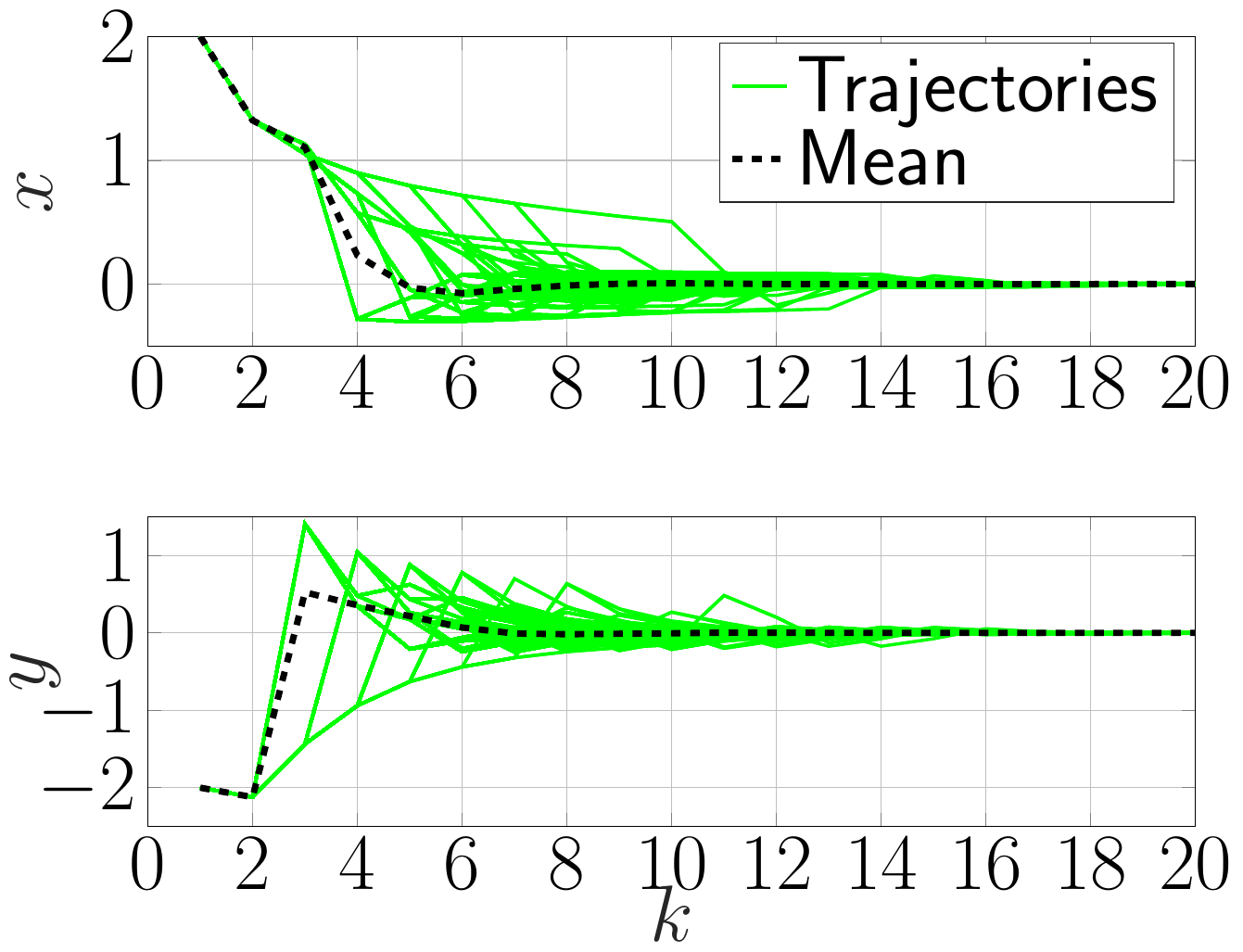}
 \includegraphics[width=0.48\columnwidth]{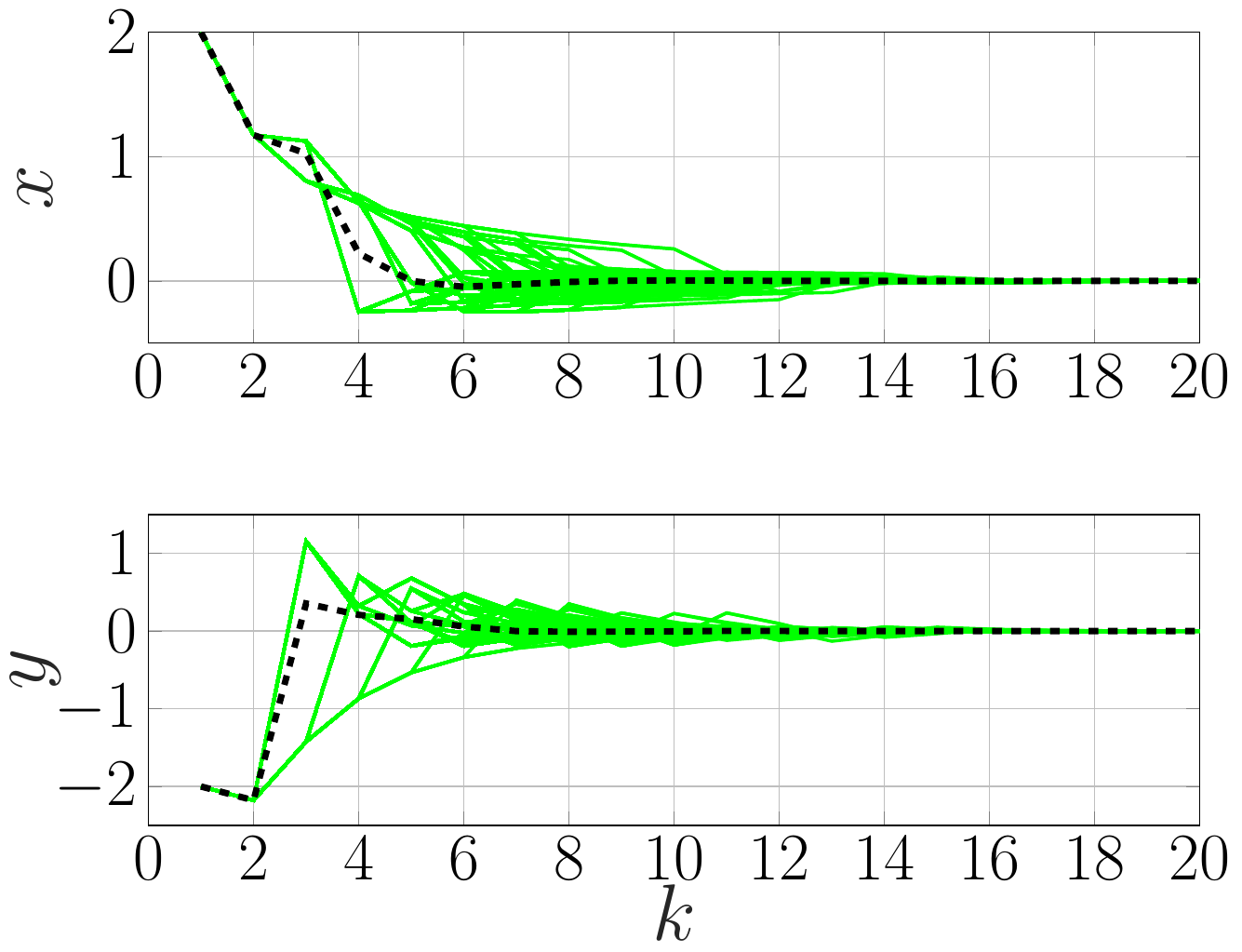}
 \caption{Trajectories of the closed-loop system 
      with risk-averse \ac{mpc} for \(N=6\) with (Left) \(\alpha=0.9\) and 
      (Right) \(\alpha=0.5\). 
      The green lines correspond to \(1000\) random 
      simulations.}%
 \label{fig:nonlinear_trajectories}
\end{figure}

Here we demonstrate the design of stabilizing risk-averse \ac{mpc} controllers 
for a nonlinear system. We consider the following nonlinear Markovian switching
system with three modes:
\begin{align}
	\smallmat{
	  x_{k+1}\\
	  y_{k+1}
	}
{}={} 
	A_{i_k}
	\smallmat{
	  x_{k}\\
	  y_{k}
	}
{}+{}
      c_{i_k}
      \smallmat{
	1 - e^{y_k}\\
	1 - e^{x_k}
      }
{}+{}
      B_{i_k}u_k.
\end{align}
The system matrices are
\begin{equation*}
	A_{1}
{}={}
	\smallmat{1 & 0.1 \\ 0.2 & 0.5},
\quad
	A_{2}
{}={}
	\smallmat{0.1 & -0.5 \\ -0.5 & 0.5}, 
\quad
	A_{3}
{}={}
	\smallmat{0.1 & -0.6 \\ 0.6 & 0.1}, 
\end{equation*}
\begin{equation*}
	B_{1}
{}={}
	\smallmat{1.6 \\ 0.6}, 
\quad 
	B_{2}
{}={}
	\smallmat{0.1 \\ 0.9},
\quad
	B_{3}
{}={}
	\smallmat{1 \\ 0},
\end{equation*}
and parameters \(c_1 {=} 0.2\), \(c_2 {=} -0.1\), \(c_3 {=} -0.3\). 
Stage-wise cost matrices are \(Q_{i} {}={}  I\) and \(R_{i} {}={} 100{}\cdot{} i\) 
for \(i {}\in{} \set{1,2,3}.\) 
The nominal and actual transition matrices are given by
\begin{equation*}
	P
{}={}
	\smallmat
	{
	  0.4 & 0.0 & 0.6 
	\\ 
	  0.6 & 0.0 & 0.4 
	\\ 
	  0.4 & 0.6 & 0.0
	},
\quad
	P'
{}={}
	\smallmat
	{
	  0.33 & 0.0 & 0.67 
	\\ 
	  0.56  & 0.0 & 0.44 
	\\  
	  0.33 & 0.67 & 0.0
	}.
\end{equation*}
The nonlinear system is constrained to be inside the box 
\(Y_1=[-2.5,2.5]^2 {}\times{} [-0.5,0.5]\) for all three modes.
Using \(m=0.5\) we compute the controller design parameters of
Thm.~\ref{thm:linearization} which are shown in Table~\ref{tbl:linearization_parameters}.
We take the terminal sets to be ellipsoidal  
\(X_{i}^f= \set{ x^{\tttop} P^f_{i} x \le r_{i} }\).
Finally, we simulate the system for different values of parameter \(\alpha\) of \(\AVARa\)
after we formulate the problem as described in Section~\ref{sec:risk-averse-oc-and-dp},
with initial condition \(x_0 = (2, -2)\) and  \(i_0 = 1\). 
Resulting system trajectories are reported in Fig.~\ref{fig:nonlinear_trajectories}. 
The proposed methodology successfully stabilizes the nonlinear system in the presence 
of uncertainty in the Markov transition matrix.  

\begin{table}[ht]
\caption{Controller design parameters}\label{tbl:linearization_parameters}
\centering
\scriptsize
\def\arraystretch{1.25}
\begin{tabular}{@{}c | c |  c  c  c @{}} 
 \multicolumn{1}{c}{}&\multicolumn{1}{c}{}&\multicolumn{3}{c}{\hspace{10pt}\(\delta_{i}\)}\\
 \(i\)  & \(\beta_{i}\) &   \(\alpha=1.0\) &    \(\alpha=0.9\) &  \(\alpha=0.5\)  \\ [0.2ex] 
 \hline 
  1 &  0.4421 &  0.2407 & 0.1783 & 0.1563 \\ 
  2 &  0.2210 &  0.3775 & 0.4121 & 0.3556 \\ 
  3 &  0.6631 &  0.1668 & 0.1130 & 0.0973 \\
\end{tabular}
\end{table}

\begin{figure}[ht]%
 \centering%
 \includegraphics[width=0.9\columnwidth]{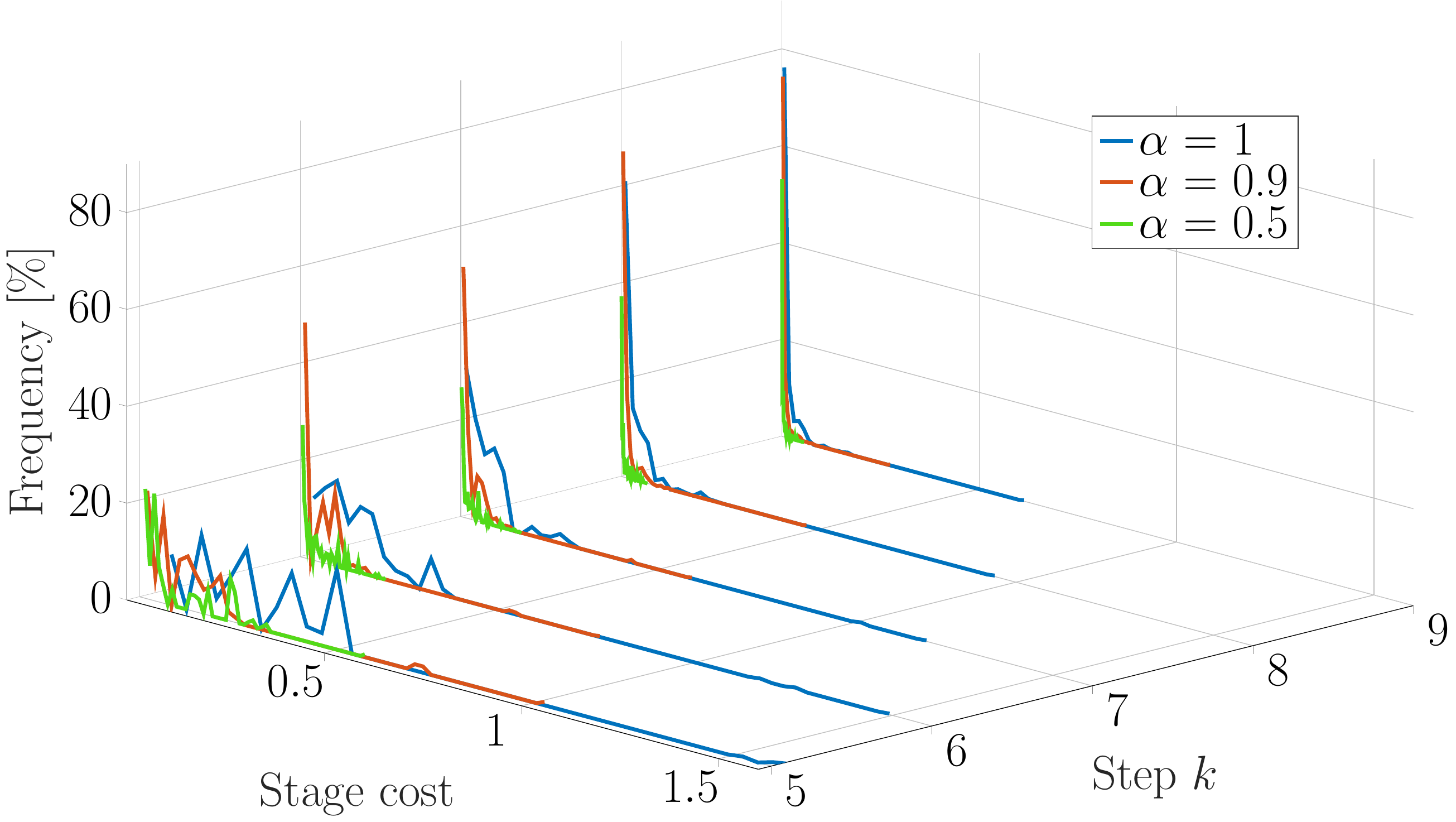}%
 \caption{Distribution of \(\ell(x_k,u_k,i_k)\) estimated using \(1000\) 
          randomly generated switching sequences. The cost of trajectories 
          corresponding to higher \(\alpha\) values are more spread out 
          compared to \(\alpha=0.5\) and have a noticeably longer right tail.}%
 \label{fig:cost_distribution}%
\end{figure}

A similar effect is observed when inspecting the distribution of 
\(\ell(x_k, u_k, i_k)\) for three \ac{mpc} controllers.
\ac{mpc} controllers with higher \(\alpha\) (closer to stochastic \ac{mpc}) 
allow for higher costs, albeit with low probability.
On the other hand, the risk-averse controller with \(\alpha = 0.5\) 
(closer to minimax \ac{mpc}) tends to produce cost distributions with shorter 
right tails. Interestingly, the point \(x_0\) is not feasible for the 
worst case controller (\(\alpha=0\)). 
The cost distributions are shown in Fig.~\ref{fig:cost_distribution}.

\section{Conclusions}
We proposed a control methodology for 
constrained nonlinear Markovian switching systems.
The proposed stability analysis framework hinges on 
\acl{dp} and leads to the formulation of 
risk-based Lyapunov-type conditions.
These conditions can be translated into an \ac{lmi} when 
the dynamics is linear, while, when the system is nonlinear
a design methodology was proposed.
In the case of \ac{mjls}, the resulting optimization problem 
can be formulated as a \ac{qcqp} and can be solved efficiently 
online enabling its use in embedded applications. 

We believe that risk-averse problems possess a favorable structure 
which can be further exploited to lead to parallelizable implementations 
akin to ones already developed for stochastic optimal control
problems~\cite{SamSopBemPat17b,SamSopBemPat17a,SamSopBemPat15}.
We plan to investigate risk-constrained 
formulations where we impose acceptable risk of 
violating the constraints instead of hard state/input constraints.
This has a potential to make the overall design
much less conservative.%
%

\begin{ack}%
This work was supported by the EU-funded H2020 project DISIRE, 
grant agreement No. 636834, the KU Leuven Research Council under
BOF/STG-15-043, the Ford-KU Leuven Research Alliance under project 
No. KUL0023 and by Research Foundation Flanders, FWO, under project
No. G086318N.
\end{ack}%

\bibliographystyle{plain}        
\bibliography{riskbib}           
\appendix
\section{Appendix}

\textit{Proof of Lemma~\ref{lem:mrses}}.
Define 
\(
	V_k 
{}\dfn{} 
	V(x_k,i_k)
\) 
and, for fixed \(x_0\in\dom V_N^\star(\cdot,i_0)\) let
\(
	x_t
{}\dfn{}
	\phi(t,x_0,i_{[t-1]})
\). 
We have%
\begin{align}\label{appendix:lem3_1}
	&\bar{\rho}_{k}
	\Bigl[V_k {-} V_0 
	{}+{} 
	c \textstyle\sum\limits_{t=0}^{k-1}\|x_t\|^2\Bigr]
{=}
	\bar{\rho}_{k}
	\Bigl[\textstyle\sum\limits_{t=0}^{k-1} 
        V_{t+1}{-}V_{t} + c \|x_{t}\|^2\Bigr] \notag\\
{}\leq{}& 
	\sum_{t=0}^{k-1} 
	\bar{\rho}_{k}
	\bigl[V_{t+1}{-}V_{t} {+} c \|x_{t}\|^2\bigr]\notag\\
{}={}&\sum_{t=0}^{k-1}
	\bar{\rho}_{t}\bigl[V_{t+1}-V_{t} + c \|x_{t}\|^2\bigr],
\end{align}%
where the inequality is because of the subadditivity of \(\rho\) (A1 and A4)
and the last equality is because \(V_{t+1}-V_{t} + c \|x_{t}\|^2\) is independent
of \(i_{t+1},\ldots,i_{k}\).
In light of Cond.~\ref{item:lyap_condition} and given that
\(
      \bar{\rho}_{t+1}[V_{t+1}-V_{t} + c \|x_{t}\|^2] 
{}={} 
      \Crho[i_0]{}
      \circ 
      \ldots
      \circ 
      \Crho
	  [i_{t}]
	  {V_{t+1}-V_{t} + c \|x_{t}\|^2}
	  []
	  [i_{t+1}]
{}={} 
      \Crho[i_0]{}
      \circ 
      \ldots
      \circ 
      \Crho
	  [i_{t}]
	  {V(f^{\kappa}(x_{t},i_{t}),i_{t+1})-V(x_{t},i_{t}) + c \|x_{t}\|^2}
	  []
	  [i_{t+1}]
{}\leq{}
      0,
\)
and because of \eqref{appendix:lem3_1} and property A2
we have that 
\(
      \bar{\rho}_k[- V_0 + c \sum_{t=0}^{k-1}\|x_t\|^2]
{}\leq{}
      \bar{\rho}_k[V_k {-} V_0 {+}
\)
\(
    c\sum_{t=0}^{k-1}\|x_t\|^2]  
{}\leq{}
   0
\).
Using properties A3 and A4, \(\bar{\rho}_k[\sum_{t=0}^{k-1}\|x_t\|^2]  \leq \nicefrac{V_0}{c}\) 
which proves the first part of Lemma~\ref{lem:mrses}.

By Cond.~\ref{item:lyap_condition},
\(
      \Crho[i_k]{V_{k+1}-V_k}[][i_{k+1}]
{}\leq{}
      -c\|x_k\|^2 
{}\leq{}
      -c\alpha_2^{-1}V_k
{}\leq{}
      -\eta V_k
\)
for some \(\eta\in (0,1)\), so
\(\Crho[i_k]{V_{k+1}}[][i_{k+1}] \leq {\beta}V_k\) with 
\(
      \beta
\dfn
      1-\eta
\in
    (0,1)
\).
We have 
\(
      \Crho[i_0]{V_1}[][i_1]
{}\leq{}
      \beta V_0
\) 
and 
\(
      \Crho[i_1]{V_2}[][i_2] 
{}\leq{}
      \beta V_1
\)%
, so 
\(
      \Crho[i_0]{
	      {\Crho[i_1]{V_2}[][i_2]}
	  }[][i_1]
{}\leq{} 
      \beta \Crho[i_0]{V_1}[][i_1]
{}\leq{} 
      \beta^2 V_0
\). 
Then,
\(
      \bar\rho_2[V_2]
{}\leq{}
      \beta^2 V_0
\)
and recursively%
\begin{align}\label{eq:pf_recursion}%
	\bar\rho_{k} \left[{}V_{k}{}\right]
{}\leq{}
	\beta^{k} V_0. 
\end{align}%
By the left hand side of Cond.~\ref{item:quad_condition}, 
\(\|x_k\|^2 \leq \nicefrac{1}{\alpha_1}V_k\) 
and applying \(\bar\rho_k\) and using~\eqref{eq:pf_recursion} and, 
subsequently the right hand side of Cond.~\ref{item:quad_condition},
\(
 \bar\rho_k(\|x_{k}\|^2) 
      \leq \bar\rho_k(\nicefrac{V_k}{\alpha_1})
      \leq \tfrac{1}{\alpha_1}\bar\rho_k(V_{k})
      \leq \tfrac{1}{\alpha_1}\beta^{k}V_0
      \leq \tfrac{\alpha_2}{\alpha_1}\beta^{k}\|x_0\|^2
\).~\hfill{\(\Box\)}

\textit{Proof of Theorem~\ref{thm:mrses-cond-tv}}.
Let \(\bar{X}\subseteq \dom V_N^\star\) be a compact \ac{ui} set.
By~\eqref{eq:dp-recursion},
\(
	V_N^\star(x,i) 
{}={}
	\Crho
	    [i]
	    {V_{N{-}1}^\star(f^{\kappa_N^\star}(x,i),j)}
	    []
	    [j]
\)
\({}+{}
	\ell(x, \kappa_N^\star(x, i), i).
\)
Then, for \((x,i) \in \bar{X}\),%
\begin{align*}
      &\Crho
	  [i]
	  {V_N^\star(f^{\kappa_N^\star}(x,i),j)}
	  []
	  [j]
{}-{}
      V_N^\star(x,i)\\
{}={}
      &\Crho
	  [i]
	  {V_N^\star(f^{\kappa_N^\star}(x,i),j)}
	  []
	  [j]
{}-{} 
      \ell(x, \kappa_N^\star(x, i), i)\\
&{}-{}
      \Crho
	    [i]
	    {V_{N{-}1}^\star(f^{\kappa_N^\star}(x,i),j)}
	    []
	    [j]\\
{}\leq{}& 
      -\ell(x, \kappa_N^\star(x, i), i)
 \leq-c\|x\|^2.
\end{align*}%
The first inequality is because  \(V_N^\star\leq V_{N-1}^{\star}\)
and property A2. 
We have that \(V_N^\star(x,i) \leq \ell_N(x,i)\leq d\|x\|^2\) for all \(x\in X^f_i\).
Because of Cond. \ref{it:mpc-mrses-Xfi-origin}, we may find \(\epsilon>0\) such that 
\(\mathcal{B}_{\epsilon} \subseteq X^f_i\),
for \(i\in\NN\). By Cond. \ref{it:mpc-mrses-VNstar-loc-bounded}, there is an \(M>d\epsilon^2\).
Then, for all \(x\in \bar{X}_i\setminus X^f_i\),  
\(
	\frac{M}{\epsilon^2}\|x\|^2
\geq
	M
\geq
	V_N^\star(x,i)
\).
Because of Cond. \ref{it:mpc-mrses-ell-bound} and the definition of \(\T\), we have that 
\(V_t^\star(x,i) \geq c \|x\|^2\) for all \((x,i)\in\dom V_t^{\star}\)
for \(t\in\N_{[1, N]}\). 
The proof is complete since \(V = V_N^\star + \delta_{\bar{X}}\) satisfies all 
conditions of Lemma~\ref{lem:mrses}.~\hfill{\(\Box\)}

\textit{Proof of Theorem~\ref{thm:linearization}}.
Define 
\(
	e(x,i) 
{}={} 
	f^{\kappa}(x,i) 
{}-{} 
	\hat{f}^\kappa(x,i)
\).
By Assumption~\ref{as:f-beta-smooth} and since 
\(
	f^\kappa(0,i){}={}0
\) 
for all \(i{}\in{}\NN\), 
\(
	\|e(x,i)\| 
{}\leq{} 
	\nicefrac{\beta_{i}}{2}\|x\|^2
\).
It is 
\(
	  \Delta(x,i) 
{}={} 
	  \Crho
	      [i]
	      {f^\kappa(x,i)^\tttop P_{j}f^\kappa(x,i)}
	      []
	      [j]
{}-{} 
	  \Crho
	      [i]
	      {x^\tttop \bar{A}_{i}^\tttop P_{j} \bar{A}_{i} x} 
	      []
	      [j]
\).
Since \(\Crho[i]{{}\cdot{}}\) is convex and monotone, it is nonexpansive with respect
to the infinity norm~\cite[p.~302]{shapirolectures}, thus for \(x\in X^f_i\)%
\begin{align*}%
&|\Delta(x,i)| 
   \leq 
    \max_{j\in\cov(i)} 
        | f^{\kappa}(x,i)^\tttop P_{j} f^{\kappa}(x,i) 
       {-} x^\tttop\bar{A}_{i}^{\tttop} P_{j} \bar{A}_{i} x
        |\\
  &\quad{}= 
    \max_{j\in\cov(i)} 
        |e(x,i)^\tttop P_j e(x,i) 
       + 2x^{\tttop} \bar{A}_{i}^{\tttop} P_{j} e(x,i)
        |\\
  &\quad{}\leq 
    \max_{j\in\cov(i)}\|P_{j}\|
        (\tfrac{\beta_{i}^2}{4}\|x\|^4 
       + \beta_{i}\|\bar{A}_{i}\| \|x\|^3
        ) \leq \sigma_i\|x\|^2
\end{align*}%
Therefore, 
\(
	  \mathcal{L}\ell_N(x,i)
{}={}
	  \mathcal{L}'\ell_N(x,i) {}+{} \Delta (x,i)
{}\leq{}
	- x^{\tttop}(\bar{Q}_{i} + (m_i - \sigma_i))x 
{}\leq{}
	-\ell(x,\kappa(x,i),i)\),
for all \(x{}\in{}X^{f}_{i}\) and since \(X^{f}\) is \ac{ui}, 
\(\T\ell_N \leq \ell_N\).
The assertion follows from Thm.~\ref{thm:mrses-cond-tv}.~\hfill{\(\Box\)}%
\end{document}